\documentclass[12pt]{article}
\usepackage{amssymb,amsfonts,amsmath, psfrag,eepic,colordvi}
\topskip  
\numberwithin{equation}{section}
\newtheorem{theorem}{Theorem}[section]

\newtheorem{lemma}[theorem]{Lemma}

\newtheorem{example}[theorem]{Example}

\begin{document}
\parskip 6pt

\pagenumbering{arabic}
\def\sof{\hfill\rule{2mm}{2mm}}
\def\ls{\leq}
\def\gs{\geq}
\def\SS{\mathcal S}
\def\qq{{\bold q}}
\def\MM{\mathcal M}
\def\TT{\mathcal T}
\def\EE{\mathcal E}
\def\lsp{\mbox{lsp}}
\def\rsp{\mbox{rsp}}
\def\pf{\noindent {\it Proof.} }
\def\mp{\mbox{pyramid}}
\def\mb{\mbox{block}}
\def\mc{\mbox{cross}}
\def\qed{\hfill \rule{4pt}{7pt}}
\def\block{\hfill \rule{5pt}{5pt}}

\begin{center}
{\Large\bf   Self-dual interval orders and row-Fishburn matrices }
\vskip 6mm
\end{center}

\begin{center}
{\small   Sherry H. F. Yan, \,\,Yuexiao Xu\\[2mm]
 Department of Mathematics, Zhejiang Normal University, Jinhua
321004, P.R. China
\\[2mm]
 huifangyan@hotmail.com
  \\[0pt]
}
\end{center}

\noindent {\bf Abstract.}   Recently, Jel\'{i}nek derived  that the
number of self-dual interval orders of reduced size $n$ is twice the
number of row-Fishburn matrices of size $n$ by using generating
functions.
  In this paper, we present a bijective proof of this
relation by establishing a  bijection  between two variations of
upper-triangular matrices  of nonnegative integers.   Using the
bijection, we provide a combinatorial proof  of the refined
relations between self-dual Fishburn matrices and
 row-Fishburn matrices in answer to a problem proposed by  Jel\'{i}nek.




\section{Introduction}
A  poset is said to be  an  {\em interval order}  ( also known as
$(2+2)$-free poset) if it does not contain an induced subposet that
is isomorphic to $2+2$, the union of two disjoint $2$-element
chains. Let $P$ be a poset with a strict order relation  $\prec$.  A
{\em strict down-set} of an element $x \in P$ is the set $D(x)$ of
all the elements of P that are smaller than $y$, i.e., $D(y) = \{y
\in  P: y \prec x\}$. Similarly, the {\em strict up-set} of $x$,
denoted by $U(x)$, is the set $\{y \in P: y\succ x\}$. A poset $P$
is $(2+2)$-free if and only if its sets of strict down-sets,
$D(P)=\{D(x): x\in P\}$ can be
 written as
  $$
  D(P)=\{D_1, D_2, \ldots, D_m \}
 $$
where $\emptyset=D_1\subset D_2\subset\ldots \subset D_m$, see
\cite{Bog, melon}. In such context, we say that $x\in P$ has {\em
level} $i$ if $D(x)=D_i$.  An element $x$ is said to be a {\em
minimal} element if $x$ has level $1$. Following Fishburn
\cite{fish1}, we call the number $m$ of distinct strict down-sets
the {\em magnitude} of P. It turns out that $m$ is also equal to the
number of distinct strict up-sets, and we can order the strict
up-sets of $P$ into a decreasing chain
$$
  D(P)=\{U_1, U_2, \ldots, U_m \}
 $$
where $ U_1\supset U_2\supset\ldots \supset U_m=\emptyset$, see
\cite{fish1, fish2}. We say that $x$ has up-level $i$ if $U(x) =
U_i$. An element $x$ is said to be a {\em maximal} element if $x$
has up-level $m$.

The {\em dual} of a poset $P$ is the poset $\overline{P}$ with the
same elements as $P$ and an order relation  $\overline{  \prec}$
defined by $x\overline{  \prec} y$ $\Longleftrightarrow$ $y\prec x$.
A poset is {\em self-dual} if it is isomorphic to its dual.

Fishburn \cite{fish1, fish3} did pioneering work on interval orders;
for instance, he showed the basic theorem that a poset is an
interval order if and only if it is $(2 + 2)$-free and  established
a bijection between interval orders and a certain kind of integer
matrices, called Fishburn matrices. Recently, Bousquet-M\'elou et
al. \cite{melon} constructed bijections between interval orders and
ascent sequences, between ascent sequences and permutations avoiding
a certain pattern, between interval orders and regular linearized
chord diagrams by Stoimenow \cite{Stoi}.  Several other papers have
focused on bijections between interval orders and other objects. For
instance, Dukes and Parviainen \cite{Duck} have described a direct
bijection between Fishburn matrices and ascent sequences, while the
papers of Claesson et al. \cite{cla1} and Dukes et al. \cite{Duck3}
extend the bijection between interval orders and Fishburn matrices
to more general combinatorial structures.

 A {\em Fishburn} matrix of
size $n$ is an upper-triangular matrix with nonnegative integers
which sum to $n$ and  each row and each column contains a nonzero
entry.   Throughout this paper that each matrix has its rows
numbered from top to bottom, and columns numbered left-to-right,
starting with row and column number one. We let $M_{i,j}$ denote the
entry of M in row $i$ and column $j$.  The {\em size} of a matrix
$M$ is the sum of all its entries. Moreover, the {\em dimension} of
an upper triangular matrix is defined to the number of rows.

The {\em dual} matrix of $M$, denoted by $\overline{M}$, is obtained
from $M$ by transposition along the diagonal running from
bottom-left to top-right. More precisely,  for $1\leq i,j\leq m$, we
have $\overline{M}_{i,j}=M_{m+1-j,m+1-i}$ where $m$ is the dimension
of $M$. If a matrix M is equal to $\overline{M}$, we call it {\em
self-dual}.

Fishburn \cite{fish1, fish3} showed that an interval order $P$ of
magnitude m corresponds to an $m\times m$ Fishburn matrix $M$ with
$M_{i,j}$ being equal to the number of elements of $P$ that have
level $i$ and up-level $j$. Jel\'{i}nek \cite{jel} showed that the
Fishburn's bijection turns out to be    a bijection between
self-dual interval orders of size $n$ and self-dual Fishburn
matrices of size $n$.

Following the terminologies given in  \cite{jel},   we distinguish
three types of cells in a Fishburn matrix $M$ of dimension $k$ : a
cell $(i, j)$ is a diagonal cell if $i+j = k+1$, i.e., $(i, j)$
belongs to the north-east diagonal of the matrix. If $i+j < k +1$
(i.e., $(i, j)$ is above and to the left of the diagonal) then $(i,
j) $ is a North-West cell, or NW-cell, while if $i + j > k + 1$,
then $(i, j)$ is an SE-cell. Clearly, NW-cells and diagonal cells
together determine a self-dual Fishburn matrix. The {\em reduced
size} of a self-dual fishburn matrix $M$ is the sum of all diagonal
cells and NW-cells. The {\em reduced size} of a self-dual interval
order $P$ is the reduced size of its corresponding self-dual
Fishburn matrix under Fishburn's bijection.

 A {\em row-Fishburn} matrix of size $n$
is defined to be an upper-triangular matrix with nonnegative
integers which sum to $n$ and each row contains a nonzero entry.
   In a matrix $A$, the {sum} of a column (resp. row) is
defined to the sum of all the entries in this column (resp. row). A
column or a row is said to be {\em zero} if it contains no nonzero
entries.   The  set of self-dual Fishburn matrices of reduced size
$n$ is denoted by $\mathcal{M}(n)$. Denote by $\mathcal{ M}(n,k,p)$
be the set of self-dual Fishburn matrices of reduced size $n$ whose
first row has sum $k$ and diagonal cells have sum $p$.    Let
$\mathcal{RM}(n)$   be the set of row-Fishburn matrices of size $n$.
The set of row-Fishburn matrices in $\mathcal{RM}(n)$ whose last
column has sum $k$ is denote by $\mathcal{RM}(n,k)$. Denote by
$\mathcal{RM}(n,k,p)$ be the set of row-Fishburn matrices in
$\mathcal{RM}(n,k)$  whose    first row has sum $p$. Moreover, the
set of self-dual interval orders of reduced size $n$ is denoted by
$\mathcal{I}(n)$.

Based on the bijection between interval orders and Fishburn
matrices, Jel\'{i}nek \cite{jel} presented a new method to derive
formulas for the generating functions of interval orders, counted
with respect to their size, magnitude, and number of minimal and
maximal elements, which generalize previous results on refined
enumeration of   interval orders obtained by Bousquet-M\'elou et al.
\cite{melon}, Kitaev and Remmel \cite{Kitaev}, and Dukes et al.
\cite{Duck2}. Applying the new method, Jel\'{i}nek \cite{jel}
obtained formulas for the generating functions of  self-dual
interval orders with respect to analogous statistics.  From the
obtained  generating functions,   relations between self-dual
Fishburn matrices  and row-Fishburn matrices were derived, that is,
\begin{equation}\label{eq1}
|\mathcal{M}(n, k, 0)|=  |\mathcal{RM}(n,k)|,
\end{equation}
 and for $p\geq 1$
\begin{equation}\label{eq2}
|\mathcal{M}(n,k,p)|=  |\mathcal{RM}(n,k,p)|.
\end{equation}
Combining   the bijection between self-dual interval orders  and
self-dual Fishburn martices,  formulas (\ref{eq1}) and (\ref{eq2}),
Jel\'{i}nek derived that for $n\geq 1$,
\begin{equation}\label{eq3}
|\mathcal{I}(n)|=|\mathcal{M}(n)|=  2|\mathcal{RM}(n)|,
\end{equation}
and asked for  bijective proofs  of  (\ref{eq1}) and (\ref{eq2}).
The main objective  of this paper is to present    bijective proofs
of these formulas by establishing a one-to-one correspondence
between two variations of upper-triangular matrices  of nonnegative
integers.

Let $\mathcal{M}(n,k)$ be the set of self-dual Fishburn matrices of
reduced size $n$ whose first row has sum $k$. Denote by
$\mathcal{EM}(n,k)$ (resp. $\mathcal{OM}(n,k)$ ) be the set of
self-dual Fishburn matrices in $\mathcal{M}(n,k)$ whose dimension
are even (resp. odd). Using the bijection between two variations of
upper-triangular matrices  of nonnegative integers, we derive that
\begin{equation}\label{eq8}
|\mathcal{EM}(n,k)|=|\mathcal{OM}(n,k)|=|\mathcal{RM}(n,k)|.
\end{equation}

\section{The bijective proofs}
Recall that a self-dual Fishburn matrix is determined by its
NW-cells and diagonal cells. Given   a self-dual Fishburn matrix
$M$,
  the {\em reduced matrix} of $M$, denoted by $R(M)$,  is a matrix
obtained from $M$ by filling all the SE-cells with zeros. An
upper-triangular matrix is said to a {\em super } triangular matrix
if all  its SE-cells are zero.

\begin{lemma}\label{le}
Let   $M'$ be a   super triangular matrix  of dimension $m$.  Then
$M'$ is a reduced matrix of a self-dual Fishburn matrix  if and only
if it satisfies the following two conditions:
\begin{itemize}
\item[(i)]for $1\leq i\leq \lceil{m\over 2}\rceil$, each column $i$ contains
a
nonzero entry;

\item[(ii)] for $1\leq i\leq \lceil{m\over 2}\rceil$,  either row $i$ or column $m+1-i$
contains a nonzero entry.
\end{itemize}
\end{lemma}
\pf  Let $M$ be a self-dual Fishburn matrix with $R(M)=M'$. Clearly,
$M'$ is a super triangular matrix. Since the first $\lceil{m\over
2}\rceil$ columns of $M'$ are the same as those in $M$, the
condition $(i)$  follows immediately. It remains to show that $M'$
satisfies condition  $(ii)$.  Since $M$ is self-dual Fishburn
matrix, for all $1\leq i\leq m$,
  row $i$ must contains  a nonzero entry,   that is, $$\sum_{j=1}^{m}
M_{i,j}=\sum_{j=1}^{m-i}
M_{i,j}+\sum_{j=m+1-i}^{m}M_{i,j}=\sum_{j=1}^{m-i}
M_{i,j}+\sum_{j=1}^{i}M_{j,m+1-i}> 0.$$  Hence, for $1\leq i\leq
\lceil{m\over 2}\rceil$,  either row $i$ or column $m+1-i$  of
$R(M)$ contains a nonzero entry.    Therefore,  the condition $(ii)$
holds for $R(M)$.

Conversely, given a   super triangular matrix $M'$  satisfying
conditions $(i)$ and $(ii)$,  We can recover a self-dual matrix $M$
from $M'$ by filling the SE-cell $(m+1-j, m+1-i)$  with $M'_{i,j}$.
      If $1\leq i\leq \lceil{m\over 2}\rceil$, the
sum of row $i$ of $M$ is given by
$$ \sum_{j=i}^{m}
M_{i,j}=\sum_{j=i}^{m-i}
M_{i,j}+\sum_{j=m+1-i}^{m}M_{i,j}=\sum_{j=1}^{m-i}
M_{i,j}+\sum_{j=1}^{i}M_{j,m+1-i}=\sum_{j=1}^{m-i}
M'_{i,j}+\sum_{j=1}^{i}M'_{j,m+1-i}. $$ By   the    condition
$(ii)$, we have $\sum_{j=i}^{m} M_{i,j}=\sum_{j=1}^{m-i}
M'_{i,j}+\sum_{j=1}^{i}M'_{j,m+1-i}>0$, which implies that row $i$
contains a nonzero entry. If $\lceil{m\over 2}\rceil+1\leq i\leq m$,
the sum of row $i$ of $M$ is given by
$$
\sum_{j=i}^{m} M_{i,j}=
\sum_{j=i}^{m}M'_{m+1-j,m+1-i}=\sum_{j=1}^{m+1-i}M'_{j,m+1-i},
$$
which implies that the sum of row $i$ of $M$ is the same as that of
column $m+1-i$ of $M'$. By condition $(i)$, row $i$ contains a
nonzero entry. Hence $M$ is a self-dual Fishburn matrix with
$R(M)=M'$. This completes the proof. \qed

 Denote by $\mathcal{SM}_k(n)$ the set of all
super triangular matrices   of size $n$  and dimension $2k+1$ having
the following two properties:
\begin{enumerate}
\item[$(a)$]    for $1\leq i\leq k$, each column $i$ contains a
nonzero entry;
 \item[$(b)$] for $1\leq i\leq k$,  either row $ k+1-i$ or
column $k+1+i$ contains a nonzero entry.
\end{enumerate}
 Let $\mathcal{SM}(n)=\bigcup_{k\geq 0}\mathcal{SM}_k(n)$.

Now we proceed to present a map $\alpha$ from $\mathcal{M}(n)$ to
$\mathcal{SM}(n)$.
 Given a  nonempty self-dual
matrix $M$ of dimension $m$,    let $\alpha(M)$  be the matrix
obtained from $M$ by the following procedure.
\begin{itemize}
\item  If $m=2k+1$ for some integer $k\geq 0$, then let $\alpha(M)$ be the
matrix obtained from the reduced matrix  $R(M)$  of $M$ by
interchanging the cell $(i,k+1)$ and the diagonal  cell $(i, m+1-i)$
for $1\leq i\leq k$.
\item If $m=2k$ for some integer $k\geq 1$, then let $A$ be the matrix
obtained from $R(M)$ by adding one zero row and one zero column
immediately after column $k$ and row $k$. Define $\alpha(M)$ to be
the matrix obtained from $A$ by interchanging the cell $(i,k+1)$ and
the diagonal  cell $(i, m+1-i)$ of the resulting matrix $A$.
\end{itemize}

Obviously,  $\alpha(M)$ is a super triangular matrix of dimension
$2k+1$ and  size $n$. It easy to check that the map $\alpha$
preserves the first $k$ columns and  the total  sum of
  row $i$ and column $m+1-i$ of the reduced matrix $R(M)$.   By Lemma \ref{le},
  the matrix $\alpha(M)$ has properties $(a)$ and $(b)$.  Hence $\alpha(M)$ is a
super triangular matrix in $\mathcal{SM}(n)$.

Conversely,  given a     super  triangular   matrix $M'$ in
$\mathcal{SM}(n)$ of dimension $2k+1$, we can recover a matrix $M\in
\mathcal{M}(n)$ with $\alpha(M)=M'$. First we   interchange  the
cell $(i,k+1)$ with the diagonal cell $(i, m-i)$ for $1\leq i\leq
k$. Then we obtain a matrix $A$ by  deleting  column $k+1$ and row
$k+1$ if they are zero.  It is easy to check that properties $(a)$
and $(b)$ ensure  that the obtained matrix $A$ is the reduced matrix
of a self-dual Fishburn matrix. Let $M$ be a self dual Fishburn
matrix with $R(M)=A$. Hence $\alpha$ is a bijection between
$\mathcal{M}(n)$ and $\mathcal{SM}(n)$.

Let $M$ be a super   triangular matrix of dimension $2k+1$, then
column $k+1$ is called a {\em center} column.  From the construction
of the bijection $\alpha$, we see that the map $\alpha$ transforms
the sum of  the diagonal cells of a self-dual matrix to the sum of
the center column of a super triangular matrix. Hence, we have the
following result.

\begin{theorem}\label{th1}
The map $\alpha$ is a bijection between $\mathcal{M}(n)$ and
$\mathcal{SM}(n)$. Moreover, the  bijection $\alpha$  preserves the
sum of the first row, and transforms the sum of  the diagonal cells
of a self-dual matrix to the sum of the center column of a super
triangular matrix.
\end{theorem}

\begin{example}
Consider a matrix $A\in \mathcal{M}(5)$,
$$
A=\begin{bmatrix}
 1 &0  &1&0  &0 \\
0 &1  &1&1  &0 \\
0 &0   &0&1  &1\\
 0 &0  &0 &1  & 0\\
  0 &0   &0&0  & 1\\
\end{bmatrix}.$$
The  reduced matrix of $A$ is given by
$$
R(A)=\begin{bmatrix}
 1 &0  &1&0  &0 \\
0 &1  &1&1  &0 \\
0 &0   &0&0  & 0\\
 0 &0  &0 &0  & 0\\
  0 &0   &0&0  & 0\\
\end{bmatrix},$$

and we have

$$
\alpha(A)=\begin{bmatrix}
1 &0  &0 &0  &1\\
0 &1  &1 &1  & 0\\
 0 &0  &0 &0  & 0\\
 0 &0  &0 &0  & 0\\
  0 &0 &0  &0  & 0\\
\end{bmatrix}.$$

\end{example}

Let $\mathcal{B}(n)$ the set of upper-triangular matrices of size
$n$ in which each row contains a nonzero entry except for the first
row.   Given  a nonempty matrix $A\in \mathcal{RM}(n)$, we can get
two distinct matrices in $\mathcal{B}(n)$  from $A$ by either doing
nothing or adding a zero row  and a zero column before the first row
and the first column.  Thus for $n\geq 1$ we have the following
relation
\begin{equation}\label{eq4}
 |\mathcal{B}(n)|=2|\mathcal{RM}(n)|.
\end{equation}
Now we proceed to  construct a bijection between the set
$\mathcal{SM}(n)$ and the set $\mathcal{B}(n)$. Before constructing
the bijection, we need some definitions. In  a matrix $A$  with $m$
rows, the operation of adding column $i$ to column  $j$ is defined
by increasing  $A_{k,j}$ by $A_{k,i}$  for each $k=1,2,\ldots, m$.

Let $\mathcal{B}(n,k,p)$ be the set of  matrices in $\mathcal{B}(n)$
whose whose first row has sum $p$ and last column has sum $k$.
Similarly,  let  $\mathcal{SM}(n,k,p)$ be the set of
 matrices in $\mathcal{SM}(n)$ whose first row has sum $k$ and
center column has sum $p$.

 \begin{theorem}\label{th2}
There is a bijection $\beta$ between  $\mathcal{SM}(n)$  and
$\mathcal{B}(n)$. Moreover, the map $\beta$ is essentially a
bijection between  $\mathcal{SM}(n,k,p)$  and $\mathcal{B}(n,k,p)$.
\end{theorem}
 \pf Given a nonempty triangular  matrix $A\in \mathcal{SM}(n) $ of
dimension $2k+1$, we recursively construct a sequence of super
triangular matrices $A^{(0)}, A^{(1)}, \ldots, A^{(l)}$.  Let
$A^{(0)}=A$ and assume that we have obtained the matrix $A^{(j)}$.
  Let $A^{(j)}$ be a super triangular matrix of dimension $2r+1$
for some integer $r$. For $1\leq i\leq r$, if each column $r+1+i$ is
zero, then let $A^{(l)}=A^{(j)}$. Otherwise, we proceed to generate
the matrix $A^{(j+1)}$ by the following {\em insertion} algorithm.
\begin{itemize}
\item  Find the largest value $i$ such that column $r+1+i$
contains a nonzero entry. Then fill the entries of column  $r+1+i$
with zeros.
\item Insert  one column immediately after column $r+1-i$, one zero
row immediately after row $r+1-i$,  one zero column immediately
before column $r+1+i$ and one zero  row  immediately before row
$r+1+i$. Let the entry in row $j$ of  the new inserted column after
column $r+1-i$ be filled with the entry in row $j$ of  column
$r+1+i$ of $A^{(j)}$ for $1\leq j\leq 2r+1$.
\end{itemize}
   Suppose that $A^{(l)}$ is
of dimension $2q+1$. Then the last $q$ rows and $q$ columns  of
$A^{(l)}$ are zero rows and columns.  Let $B$ be an upper-triangular
matrix obtained from $A^{(l)}$ by deleting the last $q$ columns and
$q$ rows.    From the    above insertion procedure to generate
$A^{(j+1)}$ form $A^{(j)}$ ,  we see that the inserted column after
column $r+1-i$ contains a nonzero entry. This ensures that each
matrix $A^{(j)}$ has property $(a)$ with $0\leq j\leq l$. Hence each
column $i$ of $A^{(l)}$ contains a nonzero entry with $1\leq i\leq
q$. Hence, $B$ is an upper-triangular matrix in which each column
contains a nonzero entry except for the last column. Moreover, the
insertion algorithm preserves the sum of each nonzero row of $A$,
which implies that $B$ is of size $n$.  Let $\beta(A)$ be the dual
matrix of $B$.  Hence we have $\beta(A)\in \mathcal{B}(n)$.

Conversely, we can construct a matrix   $A=\beta'(A')$ in
$\mathcal{SM}(n)$ from a matrix $A'$ of dimension $k+1$ in
$\mathcal{B}(n)$. Let  $B$
  be the dual matrix of $A'$. Define $M$ to be
  a matrix of dimension $2k+1$ obtained from  $B$ by adding $k$ consecutive zero rows and    $k$ consecutive
 zero columns immediately after column $k+1$ and row $k+1$.
Clearly, the obtained matrix is a super  triangular matrix having
property $(a)$. If for all $1\leq i\leq k$, either row $k+1-i$ or
column $k+1+i$ contains a nonzero entry, then we do nothing for $M$
and let $A=M$. Otherwise, we can construct a new super triangular
matrix $A$ by the following {\em removal} algorithm.
\begin{itemize}
\item Find the  least  value $i$ such that neither row $k+1-i$ nor column $k+1+i$ contains a nonzero entry.
Then we obtain a  super   triangular matrix by     adding column
$k+1-i$ to column $k+2+i$ and removing columns $k+1+i$, $k+1-i$ and
rows $k+1-i$, $k+1+i$.
 \item Repeat the above procedure for the resulting matrix until
the obtained matrix has property $(b)$.
\end{itemize}
 Obviously, the obtained matrix $A$ is a super
triangular matrix    having properties  $(a)$ and $(b)$. Since the
algorithm preserves the sums of entries in each non-zero row of $B$,
the matrix $A$ is of size $n$ and the sum of the first row of $A$ is
the same as  that of $B$.   The property  $(b)$ ensures that the
inserted columns in the insertion algorithm are the removed columns
in the removal algorithm. Thus the map $\beta'$ is the inverse of
the map $\beta$.  From the construction of the removal algorithm,
the sum of the center column of $A$ is equal to the sum of the last
column of $B$ as well as the the sum of the first row of $A'$. This
completes the proof. \qed

\begin{example}
Consider a matrix $A\in \mathcal{SM}(6)$,
$$
A=\begin{bmatrix}
1 &0  &0 &1  &1\\
0 &1  &1 &1  & 0\\
 0 &0  &0 &0  & 0\\
 0 &0  &0 &0  & 0\\
  0 &0 &0  &0  & 0\\
\end{bmatrix}.$$
Let $A^{(0)}=A$. By applying the insertion algorithm, we get

$$A^{(1)}=\begin{bmatrix}
1 &{\textbf   1}&0  &0 &1  &\textbf{0}&0\\
\textbf{0} &\textbf{0}&\textbf{0}  &\textbf{0} &\textbf{0} & \textbf{0} & \textbf{0}\\
0 &\textbf{0}&1  &1 &1  &\textbf{0}& 0\\

 0&\textbf{0} &0  &0 &0  &\textbf{0}& 0\\

 0&\textbf{0} &0  &0 &0 &\textbf{0} & 0\\
 \textbf{0}&\textbf{0} &\textbf{0}  &\textbf{0} &\textbf{0} &\textbf{0} & \textbf{0}\\
  0 &\textbf{0}&0 &0  &0 &\textbf{0} & 0\\

\end{bmatrix},$$

$$A^{(2)}=\begin{bmatrix}
1 & 1 &0 &\textbf{1} &0 &\textbf{0}&0  & 0 &0\\
 0  & 0 & 0 &\textbf{0}  & 0  &\textbf{0}& 0  &  0  &  0 \\
0 & 0 &1 &\textbf{1} &1 &\textbf{0}&0  & 0 & 0\\
\bf{0} & \textbf{0} &\textbf{0} &\textbf{0} &\textbf{0}&\textbf{0} &\textbf{0}  &\textbf{0}& \textbf{0}\\
 0& 0 &0 &\textbf{0} &0 &\textbf{0}&0  & 0 & 0\\
 \textbf{0}& \textbf{0}  &\textbf{0 }&\textbf{0} &\textbf{0} &\textbf{0}&\textbf{0} & \textbf{0}  & \textbf{0}\\

 0& 0  &0 &\textbf{0} &0 &\textbf{0}&0 & 0  & 0\\
  0 & 0  & 0 &\textbf{0}  & 0 &\textbf{0} & 0  & 0 &  0 \\
  0 & 0 &0 &\textbf{0} & 0  & \textbf{0 } &0 & 0  & 0\\

\end{bmatrix},$$
where the inserted rows and columns are illustrated in bold at each
step of the insertion algorithm. Removing the last $4$ zero rows and
$4$ zero columns, we get
$$B=\begin{bmatrix}
1 &1   &0  &1  &0       \\
0 &0   &0  &0  &0       \\
0 &0   &1  &1  &1       \\
0 &0   &0  &0  &0       \\
0 &0   &0  &0  &0       \\
\end{bmatrix}.
$$
Finally, we obtain
$$
A'= \beta(A) =\begin{bmatrix}
0 &0  &1  &0 &0       \\
0 &0   &1  &0  &1       \\
0 &0   &1  &0  &0       \\
0 &0   &0  &0  &1       \\
0 &0   &0  &0  &1      \\
\end{bmatrix}
$$
  Conversely, given $A'\in \mathcal{B}(6)$, by
applying removal algorithm, we can get $A\in \mathcal{SM}(6)$, where
the removed rows and columns are illustrated in bold at each step of
the removal algorithm.
\end{example}

Combining the bijection between self-dual interval orders and
self-dual Fishburn matrices and Theorems \ref{th1} and \ref{th2}, we
get
  a bijective  proof of $(\ref{eq3})$.

  From  Theorems \ref{th1} and \ref{th2}, we have
$$ |\mathcal{M}(n,k,p)|=|\mathcal{ SM}(n,k,p)|=|\mathcal{B}(n,k,p)|.
$$
Given a matrix $M\in \mathcal{B}(n,k,0)$, we can get a matrix $A\in
\mathcal{RM}(n,k)$ by deleting the  first row and the first column.
Conversely, given a matrix $A'\in \mathcal{RM}(n,k)$, we can obtain
a matrix $M'\in \mathcal{B}(n,k,0)$ by inserting a zero   row and a
zero column before the first row and the first column. This yields
that
\begin{equation}\label{eq5}
  |\mathcal{M}(n,k,0)|=|\mathcal{B}(n,k,0)|=|\mathcal{RM}(n,k)|.
 \end{equation}
If $p>0$, then $\mathcal{B}(n,k,p)$ is the same as
$\mathcal{RM}(n,k,p)$. Hence, if $p>0$ then we have
\begin{equation}\label{eq6}
|\mathcal{ M}(n,k,p)|=|\mathcal{B}(n,k,p)|= |\mathcal{RM}(n,k,p)|.
\end{equation}
Therefore, we get combinatorial  proofs of (\ref{eq1}) and
(\ref{eq2}), in answer to the problem posed by Jel\'{i}nek
\cite{jel}.

Now we proceed to prove (\ref{eq8}).  Given a matrix $A\in
\mathcal{EM}(n,k)$ of dimension $2m$ for some integers $m\geq1$, let
$R(A)$ be its reduced matrix. We obtain a super triangular matrix
$A'$ from $A$ by inserting a zero column and a zero row immediately
after column $m$ and row $m$. By Lemma \ref{le}, we have $A'\in
\mathcal{SM}(n,k,0)$.

Conversely, given a matrix $A'\in \mathcal{SM}(n,k,0)$ of dimension
$2m+1$ for some integer $m\geq 1$, we can recover a self-dual matrix
$A\in \mathcal{EM}(n,k)$ as follows. First, we get a super
triangular matrix $B$ from $A'$ by deleting column $m+1$ and row
$m+1$. Let $A $ be a matrix with $B=R(A )$. Obviously, we have the
matrix $A\in \mathcal{EM}(n,k)$. Hence, we get  $
|\mathcal{EM}(n,k)|=|\mathcal{SM}(n,k,0)|. $ By (\ref{eq5}), we
deduce that
\begin{equation}\label{eq7}
|\mathcal{EM}(n,k)|=|\mathcal{SM}(n,k,0)|=|\mathcal{RM}(n,k)|.
\end{equation}
 From (\ref{eq5}) and ($\ref{eq6}$), we have
$$
\begin{array}{lll}
|\mathcal{M}(n,k)|&=&|\mathcal{M}(n,k,0)|+ \sum_{p\geq
1}|\mathcal{M}(n,k,p)|\\
&=&|\mathcal{RM}(n,k)|+\sum_{p\geq
1}|\mathcal{RM}(n,k,p)|\\
&=&2|\mathcal{RM}(n,k)|.
\end{array}
$$
Meanwhile, we have
$|\mathcal{M}(n,k)|=|\mathcal{EM}(n,k)|+|\mathcal{OM}(n,k)|$. Hence,
(\ref{eq8}) follows from (\ref{eq7}).

 \noindent{\bf Acknowledgments.}   This work was supported by
the National Natural Science Foundation of China (10901141).


\begin{thebibliography}{100}
\bibitem{Bog}
K.P. Bogart, An obvious proof of Fishburn's interval order theorem,
{\em Discrete Math.} {\bf 118} (1993), 239--242.

\bibitem{melon}
 M. Bousquet-M\'elou, A. Claesson, M. Dukes, S. Kitaev,   $(2+2)$-free posets, ascent sequences and pattern avoiding permutations, {\em
J. Combin. Theory Ser. A} {\bf 117} (2010),  884--909.


\bibitem{cla1}
A. Claesson, M. Dukes, and M. Kubitzke,  Partition and composition
matrices,  {\em J.  Combin.  Theory, Ser. A}, {\bf 118} (2011),
1624--1637.
\bibitem{Duck}
M. Dukes, R. Parviainen, Ascent sequences and upper triangular
matrices containing non-negative integers, {\em Electronic  J.
combin. } {\bf 17 } (2010),  R53.

\bibitem{Duck2}
M. Dukes, S. Kitaev, J. Remmel, and E. Steingr\'imsson, Enumerating
(2+2)-free posets by indistinguishable elements, arXiv:1006.2696,
2010.
\bibitem{Duck3}
M. Dukes, V. Jel\'inek, and M. Kubitzke,  Composition matrices,
(2+2)-free posets and their specializations, {\em Electronic J.
Combin.}, {\bf 18}  (2011), P44.


\bibitem{fish1}
P. C. Fishburn,    Interval lengths for interval orders: A
minimization problem,   {\em Discrete Mathematics}, {\bf 47} (1983),
63--82.
\bibitem{fish2}
P. C. Fishburn,  Interval graphs and interval orders,  {\em Discrete
Mathematics}, {\bf 55} (1985),  135--149.
\bibitem{fish3}
P. C. Fishburn, {\em Interval orders and interval graphs: A study of
partially ordered sets},  John Wiley \& Sons, 1985.

\bibitem{jel}
V.  Jel\'{i}nek, Counting self-dual interval orders,  
arXiv:1106.2261, 2011.



\bibitem{Kitaev}
S. Kitaev,  J. Remmel,  Enumerating $(2+2)$-free posets by the
number of minimal elements and other statistics,  {\em Disctere
Appl. Math.}, {\bf 159 } (2011), 2098--2108.
\bibitem{Stoi}
A. Stoimenow, Enumumeration of chord diagrams and an upper bound for
Vassiliev invariants, {\em J. Knot Theory Ramifications}  {\bf 7}
(1998), 93--114.
\bibitem{zag}
D. Zagier, Vassiliev invariants and a stange identity related to the
Dedeking eta-function, {\em Topology}   {\bf 40} (2001), 945--960.
\end{thebibliography}
\end{document}